\documentclass{amsart}[12pt]
\usepackage{amsmath, amsthm, amscd, amsfonts,}
\usepackage{graphicx,float}
\usepackage[utf8]{inputenc}
\usepackage{amssymb}

\usepackage{pb-diagram}

\newcommand{\PP}{{\mathbb{P}}}

\def\MQB{{\mathbb{Q}}}

\newtheorem{theorem}{Theorem}[section]
\newtheorem{lemma}[theorem]{Lemma}

\newtheorem{remark}[theorem]{Remark}
\newtheorem{question}[theorem]{Question}

\newtheorem{claim}[theorem]{Claim}
\numberwithin{equation}{section}

\def\MQB{{\mathbb{Q}}}

\def\rmark{\mbox{$\rm\bf\rule{0.06em}{1.45ex}\kern-0.05em R$}}
\def\pmark{\mbox{$\rm\bf\rule{0.06em}{1.45ex}\kern-0.05em P$}}
\def\nmark{\mbox{$\rm\bf\rule{0.06em}{1.45ex}\kern-0.05em N$}}
\def\vdash{\mbox{$\rm\| \kern-0.13em -$}}
\newcommand{\lusim}[1]{\smash{\underset{\raisebox{1.2pt}[0cm][0cm]{$\sim$}}
{{#1}}}}

\usepackage{pb-diagram}

\def\rmark{\mbox{$\rm\bf\rule{0.06em}{1.45ex}\kern-0.05em R$}}
\def\pmark{\mbox{$\rm\bf\rule{0.06em}{1.45ex}\kern-0.05em P$}}
\def\nmark{\mbox{$\rm\bf\rule{0.06em}{1.45ex}\kern-0.05em N$}}
\def\vdash{\mbox{$\rm\| \kern-0.13em -$}}

\begin{document}

\title[Cardinal collapsing and product forcing]{Cardinal collapsing and product forcing}

\author[M. Golshani and R. Mohammadpour]{Mohammad Golshani and Rahman Mohammadpour}

\thanks{The first author's research has been supported by a grant from IPM (No. 1400030417).
}

\thanks{The authors thank Radek Honzik for some useful comments and corrections in an earlier version of the paper.}
 \maketitle

\begin{abstract}
Suppose $\kappa$ is a singular strong limit cardinal of countable cofinality
and let $\langle \kappa_{n}: n<\omega \rangle$  be an incrasing sequence of regular cardinals cofinal in $\kappa$.
We show that if $cf(2^\kappa)= \kappa^+$, then forcing with the full product $\prod_{n<\omega}Add(\kappa_n,1)$  collapses $2^\kappa$ into $\kappa^+$. This result gives a consistent positive answer to a question of Sy Friedman. We also give a new proof of a result due to Shelah by showing that  if the sequence carries a scale of length $\kappa^+,$
then forcing with  $\prod_{n<\omega}Add(\kappa_n,1)$  adds a generic filter for $Add(\kappa^+, 1)$, and indeed
\[
\prod_{n<\omega}Add(\kappa_n,1)/fin \simeq Add(\kappa^+, 1).
\]
\end{abstract}

\maketitle
\section{Introduction}
Suppose $\langle \kappa_{n}: n<\omega \rangle$ is an increasing sequence of regular cardinals cofinal in $\kappa$.
In \cite{honzik}, Sy Friedman and Radek Honzik observed that if
 $\prod_{n<\omega}\kappa_n$ carries a scale of length $\kappa^+,$ then $\prod_{n<\omega} Add(\kappa_n, 1)$ collapses $2^\kappa$ into $\kappa^+.$

On the other hand, answering a question of Friedman and Rene David, Saharon Shelah \cite{shelah} showed that
if $\prod_{n<\omega}\kappa_n$ carries a scale of length $\kappa^+,$
then forcing with $\prod_{n<\omega} Add(\kappa_n, 1)$ adds a generic for $Add(\kappa^+, 1)$ over $V$.

As forcing with  $Add(\kappa^+, 1)$ collapses $2^\kappa$ into $\kappa^+,$ the Friedman-Honzik's result follows from
Shelah's theorem. In proofs of both  results, the assumption that $\prod_{n<\omega}\kappa_n$ carries a scale of length $\kappa^+$ seems to be essential. In personal communication \cite{friedman}, Sy Friedman asked the first author if we can remove the assumption of the existence of scale from his result with Honzik. More precisely, he asked if the following is true:
\begin{question}
Suppose $\langle \kappa_{n}: n<\omega \rangle$ is an increasing sequence of inaccessible cardinals cofinal in $\kappa.$ Does forcing with $\prod_{n<\omega} Add(\kappa_n, 1)$ collapse $2^\kappa$ into $\kappa^+?$
\end{question}
We give a consistent positive answer to Friedman's question, by proving the following:
\begin{theorem}
\label{main theorem1}
Assume  $\kappa$ is a singular strong limit cardinal of countable cofinality and  $cf(2^\kappa)=\kappa^+$.
Let $\langle \kappa_{n}: n<\omega \rangle$ be any increasing sequence of regular cardinals cofinal in $\kappa$ and let $\langle \mathbb{P}_n: n<\omega \rangle$ be a sequence of non-trivial separative forcing notions, such that each $\mathbb{P}_n$ is $\kappa_n$-closed and of size $< \kappa.$ Further suppose that for each $\lambda < \kappa_n$, each decreasing sequence $\langle p_\alpha: \alpha < \lambda \rangle$ in $\mathbb{P}_n$
 has a greatest lower bound in $\mathbb{P}_n$. Then $\prod_{n<\omega}\mathbb{P}_n$  collapses $2^\kappa$ into $\kappa^+.$
\end{theorem}
We also give a new proof of the above result of Shelah, indeed we prove the following.
\begin{theorem}
\label{main theorem2}
Assume  $\kappa$ is a singular strong limit cardinal of countable cofinality and
let $\langle \kappa_{n}: n<\omega \rangle$ be any increasing sequence of regular cardinals cofinal in $\kappa$ which carries a scale of length $\kappa^+$.  Then $\prod_{n<\omega}Add(\kappa_n,1)/fin \simeq Add(\kappa^+, 1)
$, in particular
 forcing with $\prod_{n<\omega}Add(\kappa_n, 1)$ adds a generic filter for $Add(\kappa^+, 1)$.
\end{theorem}

The structure of the paper is as follows. In Section \ref{Proof of Theorem 1.2}, we present a proof of  Theorem \ref{main theorem1}
and in Section \ref{maintheorem2}, we prove
 Theorem \ref{main theorem2}.

\section{Proof of Theorem \ref{main theorem1}}
\label{Proof of Theorem 1.2}
In this section, we give a proof of Theorem \ref{main theorem1}.  Let $\mathbb{P}=\prod_{n<\omega}\mathbb{P}_n,$ and let  $\mathbb{Q}=(\prod_{n<\omega}\mathbb{P}_n/fin).$ Also let $\pi: \mathbb{P} \rightarrow \mathbb{Q}$ be defined in the natural way:
\begin{center}
$\pi(\langle p_n: n<\omega   \rangle)=[\langle p_n: n<\omega   \rangle]/fin,$
\end{center}
where $[\langle p_n: n<\omega   \rangle]/fin$ denotes the equivalence class of $\langle p_n: n<\omega   \rangle$ in $\mathbb{Q}.$ The following can be proved easily.
\begin{lemma}
\label{projection}
$\pi$ is a projection, i.e.,
\begin{enumerate}
\item
$\pi(1_{\mathbb{P}})=1_{\mathbb{Q}}$.
\item
$\pi$ is order preserving.
\item
If $[p]/fin\leq_{\mathbb{Q}} [q]/fin$, then there exists $r \leq_{\mathbb{P}} q$ such that $[r]/fin\leq_{\mathbb{Q}} [p]/fin$.
\end{enumerate}
\end{lemma}
We now show that forcing with $\mathbb{Q}$ adds a new $\kappa^+$-sequence of ordinals. We need the following lemma.
\begin{lemma} (\cite{golshani-shelah})
\label{golshelahthm}
Assume $cf(2^\kappa)=\kappa^+$  and $\MQB$ is a $(\kappa+1)$-strategically closed forcing notion of size $2^\kappa$, such that player
II has a winning strategy where at limit stages he chooses the greatest lower bound of the previously chosen sequence.
 Then
forcing with $\MQB$ adds a new sequence of ordinals of length $\kappa^+.$
\end{lemma}
It is not difficult to show that the forcing notion $\mathbb{Q}$ is $(\kappa+1)$-strategically closed and there exists a
winning strategy for player II where at limit stages, he chooses the greatest lower bound of the previously chosen sequence.
 It follows from Lemma \ref{golshelahthm} that forcing with $\mathbb{Q}$ adds a new $\kappa^+$-sequence of ordinals,
 as requested. To complete the proof we
also need the following result.
\begin{lemma}(\cite{golshani-hayut})
\label{golhayut}
Let $\mathbb{Q}$ be a $(\kappa+1)$-strategically closed forcing notion of size $2^{\kappa}$. Let $o(\mathbb{Q})$ be the least cardinal $\mu$, such that forcing with $\mathbb{Q}$ adds a new $\mu$-sequence of ordinals (or equivalently of elements of $V$). Then forcing with $\mathbb{Q}$ collapses $2^{\kappa}$ onto $o(\mathbb{Q})$.
\end{lemma}
\begin{remark}
In \cite{golshani-hayut}, the lemma is not stated as above, but the proof and remarks after it show that the above stronger result holds.
\end{remark}
We are now ready to complete the proof of Theorem \ref{main theorem1}. By
the above lemmas, forcing with $\mathbb{Q}$ collapses $2^\kappa$
to $\kappa^+$  and by
Lemma \ref{projection}, $V^{\mathbb{Q}} \subseteq V^{\mathbb{P}},$ hence forcing with $\mathbb{P}$ collapses $2^\kappa$
to $\kappa^+$ as well.

\section{Proof of Theorem \ref{main theorem2}}
\label{maintheorem2}
In this section we prove Theorem \ref{main theorem2}.
The proof is given in two stages. At the first stage we show that forcing with $\prod_{n<\omega}Add(\kappa_n, 1)$
collapses $2^\kappa$ into $\kappa^+.$ In the next stage,
we analyze the forcing notion $\prod_{n<\omega}Add(\kappa_n, 1)/fin$
and use our results to conclude the theorem.

{\bf Stage 1:}
In this stage we show that forcing with $\prod_{n<\omega}Add(\kappa_n, 1)$
collapses $2^\kappa$ into $\kappa^+$.
Let $\vec{f}=\langle f_\alpha: \alpha<\kappa^+ \rangle$ be a scale in $\prod_{n<\omega}\kappa_n.$

Let $\mathcal{F}=\{f\in \prod_{n<\omega}h(n): f=^*f_\alpha, $ for some $\alpha<\kappa^+\}.$ Then $|\mathcal{F}|=\kappa^+,$ and it is cofinal in $(\prod_{n<\omega}\kappa_n, \leq ).$ Also let $G_n: \kappa_n \rightarrow 2$ be the Cohen generic function, added by $Add(\kappa_n, 1).$
For each $f\in \mathcal{F},$ define $g_f: \kappa \rightarrow 2,$ so that for each  $n < \omega,$
$$g_f(\kappa_{n-1}+\xi)=G_{n+1}(f(n+1)+\xi),$$ where $\kappa_{n-1} \leq \xi<\kappa_{n},$ and $\kappa_{-1}=0.$
We show that for each $g:\kappa \rightarrow 2, g\in V,$ there is $f\in \mathcal{F},$ such that $g=g_f.$

Indeed let $g\in V$ be as above and set
$$D=\{p \in \prod_{n<\omega}Add(\kappa_n, 1): \exists f \in\mathcal{F}, \forall n, \forall \xi< \kappa_n [g(\xi)=p(n+1)(f(n+1)+\xi) ]              \}.$$
Then $D$ is dense in $\prod_{n<\omega}Add(\kappa_n, 1)$. To see this, let $p \in \prod_{n<\omega}Add(\kappa_n, 1)$. By extending $p$, we
may assume that for each $n<\omega, p(n): \zeta_n \to 2$, for some $\zeta_n<\kappa_n$.
It then follows that
$\langle \zeta_n: n<\omega \rangle  \in \prod_{n<\omega}\kappa_n$.
Pick $\alpha < \kappa^+$ such that $g <^* f_\alpha$. It follows that $g < f$ for some $f \in \mathcal{F}.$
Now define $q \in \prod_{n<\omega}Add(n, 1)$ by $q(n): f(n) \to 2$, $q(n) \supseteq p(n)$
and for all $\xi < h(n), q(n+1)(f(n+1)+\xi)=g(\xi).$
Then $q$ is well-define, it extends $p$ and $q \in D$.

It follows that for some $f\in \mathcal{F},$ $g=g_f,$
as requested.

{\bf Stage 2:} In this stage we complete the proof of Theorem \ref{main theorem2}. For each $n<\omega$  set $\mathbb{P}_n=Add(\kappa_n, 1)$.  Let $\mathbb{P}=\prod_{n<\omega}\mathbb{P}_n$  and $\mathbb{Q}=\prod_{n<\omega}\mathbb{P}_n/fin.$ Let also $\pi: \mathbb{P} \rightarrow \mathbb{Q}$ be defined as in the previous section. The next lemma can be proved easily.
\begin{lemma}
\label{quotientforcingproperties}
\begin{enumerate}
\item[(a)] $\mathbb{Q}$ is $< \kappa^+$-strategically closed.
\item[(b)] the quotient forcing $\mathbb{P} / \dot{G}_{\mathbb{Q}}$ is $\kappa^+$-c.c.
\end{enumerate}
\end{lemma}
As clearly forcing with $\mathbb{P}$ preserves cardinals $\leq \kappa^+$
and by Stage 1 it collapses $2^\kappa$ into $\kappa^+$, it follows from Lemma
\ref{quotientforcingproperties}(b) that it is forcing with $\mathbb{Q}$
that collapse $2^\kappa$ onto $\kappa^+$.
Now we need the following known result:
\begin{lemma} (see \cite{jech})
\label{colforcing}
Suppose $\kappa < \lambda$ are infinite cardinals and  $\lambda^{\kappa}=\lambda$. Suppose $\mathbb{Q}$ is a $<\kappa^+$-strategically closed forcing notion of size $\lambda$ and suppose that forcing with $\mathbb{Q}$ collapses $\lambda$ into $\kappa^+$. Then
$\mathbb{Q} \simeq Col(\kappa^+, \lambda).$
\end{lemma}
By lemmas \ref{colforcing} and \ref{quotientforcingproperties}(a),
\[
\mathbb{Q} \simeq Add(\kappa^+, 1).
\]
Now by Lemma \ref{projection}, forcing with $\mathbb{P}$
adds a generic for $\mathbb{Q}$, which completes the proof of Theorem \ref{main theorem2}.

\bigskip
School of Mathematics, Institute for Research in Fundamental Sciences (IPM), P.O. Box:
19395-5746, Tehran-Iran.

E-mail address: golshani.m@gmail.com\\

Institut für Diskrete Mathematik und Geometrie, TU Wien,
1040 Vienna, Austria.

E-mail address: rahmanmohammadpour@gmail.com

\end{document}